\documentclass[11pt]{amsart}
\usepackage{amsmath,amsthm,epsfig,amssymb}
\usepackage{xypic}
\input xy
\title{Ehrhart polynomial for lattice squares, cubes and hypercubes}
\author{Eugen J. Ionascu}
\curraddr{Department of Mathematics\\ Columbus State University\\4225 University Avenue\\
Columbus, GA 31907\\
Honorific Member of the Romanian Institute of Mathematics ``Simion
Stoilow" } \email{ionascu@columbusstate.edu, math@ejionascu.ro}   \subjclass{52C07, 05A15, 68R05}
\date{March $17^{th}$, 2016}
  \keywords{Ehrhart polynomial, linear Diophantine
equations, lattice square, lattice cube, quaternions, icube, Pythagorean quadruple, twin vectors}
\begin{document}
\def\sms{\small\scshape}
\baselineskip18pt
\newtheorem{theorem}{\hspace{\parindent}
T{\scriptsize HEOREM}}[section]
\newtheorem{proposition}[theorem]
{\hspace{\parindent }P{\scriptsize ROPOSITION}}
\newtheorem{corollary}[theorem]
{\hspace{\parindent }C{\scriptsize OROLLARY}}
\newtheorem{lemma}[theorem]
{\hspace{\parindent }L{\scriptsize EMMA}}
\newtheorem{definition}[theorem]
{\hspace{\parindent }D{\scriptsize EFINITION}}
\newtheorem{problem}[theorem]
{\hspace{\parindent }P{\scriptsize ROBLEM}}
\newtheorem{conjecture}[theorem]
{\hspace{\parindent }C{\scriptsize ONJECTURE}}
\newtheorem{example}[theorem]
{\hspace{\parindent }E{\scriptsize XAMPLE}}
\newtheorem{remark}[theorem]
{\hspace{\parindent }R{\scriptsize EMARK}}
\renewcommand{\thetheorem}{\arabic{section}.\arabic{theorem}}
\renewcommand{\theenumi}{(\roman{enumi})}
\renewcommand{\labelenumi}{\theenumi}
\newcommand{\Q}{{\mathbb Q}}
\newcommand{\Z}{{\mathbb Z}}
\newcommand{\N}{{\mathbb N}}
\newcommand{\C}{{\mathbb C}}
\newcommand{\R}{{\mathbb R}}
\newcommand{\F}{{\mathbb F}}
\newcommand{\K}{{\mathbb K}}
\newcommand{\D}{{\mathbb D}}
\def\phi{\varphi}
\def\ra{\rightarrow}
\def\sd{\bigtriangledown}
\def\ac{\mathaccent94}
\def\wi{\sim}
\def\wt{\widetilde}
\def\bb#1{{\Bbb#1}}
\def\bs{\backslash}
\def\cal{\mathcal}
\def\ca#1{{\cal#1}}
\def\Bbb#1{\bf#1}
\def\blacksquare{{\ \vrule height7pt width7pt depth0pt}}
\def\bsq{\blacksquare}
\def\proof{\hspace{\parindent}{P{\scriptsize ROOF}}}
\def\pofthe{P{\scriptsize ROOF OF}
T{\scriptsize HEOREM}\  }
\def\pofle{\hspace{\parindent}P{\scriptsize ROOF OF}
L{\scriptsize EMMA}\  }
\def\pofcor{\hspace{\parindent}P{\scriptsize ROOF OF}
C{\scriptsize ROLLARY}\  }
\def\pofpro{\hspace{\parindent}P{\scriptsize ROOF OF}
P{\scriptsize ROPOSITION}\  }
\def\n{\noindent}
\def\wh{\widehat}
\def\eproof{$\hfill\bsq$\par}
\def\ds{\displaystyle}
\def\du{\overset{\text {\bf .}}{\cup}}
\def\Du{\overset{\text {\bf .}}{\bigcup}}
\def\b{$\blacklozenge$}

\def\eqtr{{\cal E}{\cal T}(\Z) }
\def\eproofi{\bsq}

\begin{abstract} In this paper we are constructing integer lattice squares, cubes or hypercubes in $\mathbb R^d$ with $d\in \{2,3,4\}$. We find
a complete description of their Ehrhart polynomial.  We characterize all the integer squares in $\mathbb R^4$, in terms of two  Pythagorean quadruple representations
of the form $a^2+b^2+c^2=d^2$, and then prove a parametrization in terms of two quaternions of all such squares. We introduce the sequence of {\em almost perfect squares in dimension n}. In dimension two, this is very close to the sequence A194154 (in OEIS). 
\end{abstract} \maketitle
 
\section{INTRODUCTION}
$\rm Eug\grave{e}ne$ Ehrhart (\cite{ee},\cite{ee2}) proved that given a $d$-dimensional compact simplicial complex in $\mathbb
R^n$ ($1\le d\le n$), denoted here generically by $\cal P$,  whose vertices are in the lattice $\mathbb Z^n$,
there exists a polynomial  $L({\cal P}, t)\in \mathbb
\mathbb Q [t]$ of degree $d$, associated with $\cal P$, satisfying

$$L({\cal P},t)=\text{the \ cardinality\  of} \ \{t\cal P\}\cap \mathbb Z^n, \ t\in \mathbb N.$$

It is very interesting that one can say more about three of the coefficients of $L({\cal P},t)$:

\begin{equation} \label{ehrhartpolynomial}L({\cal P}, t) = Vol ({\cal P}) t^d + \frac{1}{2}Vol (\partial {\cal P})t ^{d-1} +... + \chi ({\cal P}),\end{equation}

\n where  $Vol ({\cal P})$ is the usual volume of $\cal P$ normalized
with respect to the sublattice containing ${\cal P}$,  $Vol (\partial {\cal P})$  is the surface area of $\cal P$ normalized
with respect to the sublattice on each face of $\cal P$ and $\chi({\cal P})$ is the Euler characteristic of $\cal P$ (in the sense of polytopal complexes, so for a convex polytope it is equal to one). In general, the other coefficients
in (\ref{ehrhartpolynomial}) may have complicated expressions in terms of the vertices of  $\cal P$.

It is also known that the number of points in the interior of $t \cal P$ is given by $(-1)^d L({\cal P}, -t)$.
For a polytope that is a cross product of two polytopes, the Ehrhart polynomial is the product of the corresponding smaller degree Ehrhart polynomials.
In this article, we are studying this polynomial for lattice squares, cubes, and hypercubes in $\mathbb R^n$.
Such objects have been constructed and counted in several works (see \cite{gkms}, \cite{spira}, \cite{Hurlimann}, \cite{kissKutas}).

\section{Squares}

\subsection{Squares in two dimensions}

How general can a lattice square in two dimensions look like? We may assume that one of the vertices is at the origin.
Then the other vertices, in counterclockwise order, are given by $A(a,b)$, $B(a-b,a+b)$ and $C(-b,a)$ where $a$ and $b$ are integers, not both zero, and such that
$\gcd(a,b)=1$. Using the facts mentioned in the Introduction,  the polynomial in (\ref{ehrhartpolynomial}) is simply $E_\Box(t)=(a^2+b^2)t^2+2t+1$ and of course,
$E_{\overset{\circ}{\Box}}(t)=(a^2+b^2)t^2-2t+1$, $t\in \mathbb N$.

\begin{figure}
\centering
\includegraphics[scale=.3]{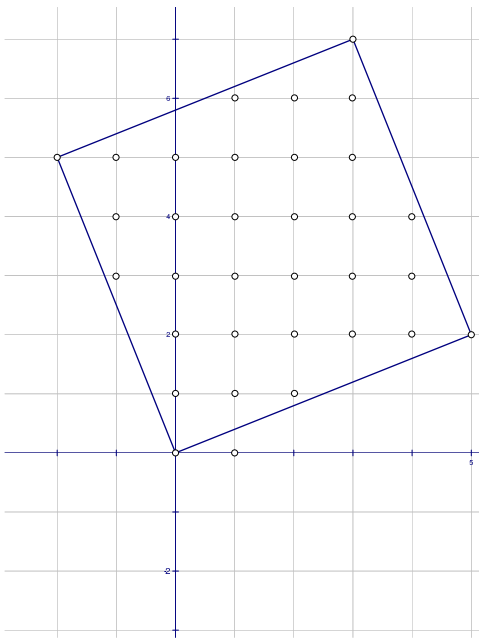}
\caption{Square with a=5 and b=2}
\label{fig:fig1}
\end{figure}

In Figure~\ref{fig:fig1}, we see the $28$ points inside of the square $OABC$, with $A(5,2)$, $B(3,8)$ and $C(-2,5)$.
One interesting problem here is to determine the sequence of possible lattice points in the interior of such a square.
The first hundred terms of this sequence, listed in increasing order, are
included in the table below:

$${\tiny
\begin{tabular}{|c|c|c|c|c|c|c|c|c|c|c|c|c|c|c|c|c|c|c|c|}
  \hline
   0 & 1 & 4 & 5 & 9 & 12 & 13& 16 & 17 & 24 & 25 & 28 & 33 & 36 & 37 & 40 & 41 & 49 & 52 & 57  \\
  \hline 60 & 61 & 64 & 65 & 72 & 73 & 81 & 84 & 85 & 88 & 96 & 97 & 100 & 101 & 105 & 108 & 112 & 113 & 116 & 121  \\
  \hline 124& 129  & 133 & 136 & 144 & 145 & 148 & 153 & 156 & 161 & 168 & 169 & 172 & 177 & 180 & 181 & 184 & 192 & 193 & 196  \\
  \hline 197& 201 & 204 & 209 & 217 & 220 & 221 & 225 & 228 & 229 & 232 & 240 & 241 & 249 & 256 & 257 & 264 & 265 & 268 & 273  \\
  \hline 276& 280 & 288 & 289 & 292 & 293 & 297 & 301 & 304 & 305 & 312 & 313 & 316 & 324 & 325 & 328 & 336 & 337 & 345 & 348  \\
  \hline
\end{tabular}}$$

\begin{figure}
\centering
\includegraphics[scale=.3]{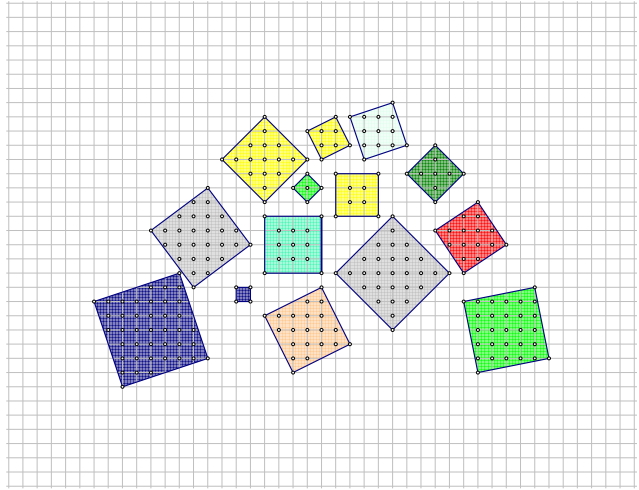}
\caption{Various squares}
\label{fig:fig2}
\end{figure}

Let us call this sequence the {\em almost perfect squares} sequence.
In Figure~\ref{fig:fig2} we see some of the squares that show that the above numbers are indeed in the sequence.
As a curiosity, $2015$ is not in this sequence but $2016$ is. This sequence is very close to A194154 in the
OEIS (On-Line Encyclopedia of Integer Sequences), but $20$ is the first term that is not in our sequence.

Let us give a reformulation of the exact form of the Ehrhart polynomial above in conjunction with the reciprocity property (about interior points), which can by shown independently using Pick's Theorem.

\begin{proposition}\label{firstpropverytrivial}  Given $a$, $b$ with $\gcd(a,b)=1$ and $t\in \mathbb N$, then
$$\sharp\{(x,y)\in \mathbb Z^2:\  ax+by, \ ay-bx \in [1,t(a^2+b^2)-1] \}=(a^2+b^2)t^2-2t+1.$$
\end{proposition}

\proof. \ Counting the points inside the square $OABC$, with $O(0,0)$, $A(a,b)$, $B(a-b,a+b)$ and $C(-b,a)$, is equivalent to count the solutions $(x,y)$ of the system
$(x,y)=\alpha (ta,tb)+\beta (-tb,ta)$ with $\alpha,\beta\in (0,1)$. Since $\alpha=(x,y)\cdot (a,b)/{t(a^2+b^2)}=\frac{ax+by}{t(a^2+b^2)}$
and $\beta=(x,y)\cdot (-b,a)/{t(a^2+b^2)}=\frac{ay-bx}{t(a^2+b^2)}$ we see that the constraints on $\alpha$ and $\beta$ is equivalent to $ax+by, \ ay-bx \in [1,t(a^2+b^2)-1]$. Hence, the result follows from the Ehrhart polynomial expression.\eproof

So, every term of the sequence of almost perfect squares is the answer to a counting as in Proposition~\ref{firstpropverytrivial}, which in particular gives the justification of one of our earlier claims:

$$\sharp\{(x,y)\in \mathbb Z^2:\  44x+9y, \ 44y-9x \in [1,2016] \}=2016.$$

\subsection{Squares in three dimensions}

In $\mathbb R^3$, we can obtain a lattice square by taking two orthogonal vectors with integer coordinates of the same length:
$u=(a,b,c)$, $v=(a',b',c')$ such that $a^2+b^2+c^2=a'^2+b'^2+c'^2=\ell $ and $aa'+bb'+cc'=0$. In \cite{gkms}, two such vectors are call {\em twin} vectors.
The number of such twin vectors having a given length is calculated in \cite{gkms}. In \cite{kissKutas}, this concept is generalized to $m$ dimensions
and called $m$-icube (a set of $m$ vectors in $\mathbb Z^n$, mutually orthogonal and of the same length). Since we are interested in generalizing these ideas to any dimension, we are going to use this terminology also.

Although it does not make much of a difference it is natural to assume that the square is {\em irreducible}, i.e., $\gcd(a,b,c,a',b',c')=1$. If we set $d=\gcd(a,b,c)$,  $d'=gcd(a',b',c')$ and $D=\gcd(bc'-b'c,ac'-a'c,ab'-b'a)$, we have the following more general formula, than in two dimensions, for the  $E_\Box(t)$.

\begin{theorem}\label{epford3} The Ehrhart polynomial of a lattice square embedded into  $\mathbb R^3$, described above and with the notation introduced is given by
\begin{equation}\label{ehrhart1}
E_\Box(t)=Dt^2+(d+d')t+1.
\end{equation}
\end{theorem}

\proof. \ The fundamental domain of the lattice containing the triangle has a volume equal to $\sqrt{n_1^2+n_2^2+n_3^2}$ where
$\overrightarrow{n}=[n_1,n_2,n_3]$ is a vector normal to the plane containing the triangle and such that $\gcd(n_1,n_2,n_3)=1$ (see \cite{abarvinok1} and \cite{ejilps}). Such a vector is given clearly by the
cross-product of $u$ and $v$: $n=\frac{1}{D}u\times v=\frac{1}{D}(bc'-b'c,-(ac'-a'c),ab'-b'a)$. Let us observe that because of Lagrange's Identity
$\ell^2=(a^2+b^2+c^2)(a'^2+b'^2+c'^2)=(aa'+bb'+cc')^2+(bc'-b'c)^2+(ac'-a'c)^2+(ab'-b'a)^2$, we conclude that $|n|=\frac{\ell}{D}$.
Hence, we see that the first coefficient of the  Ehrhart polynomial is $\sqrt{\ell}^2/(\frac{\ell}{D})=D$. Similar argument goes for the
second coefficient.\eproof

\n {\bf Standard Examples:} Let us analyze some examples which are not  $xy$, $xz$ or $yz$-plane examples. First, let us take $u=(3,-3,0)$ and $v=(1,1,4)$; observe that
$|u|^2=|v|^2=18$ and $u\cdot v=0$.  In this case $n=6(-2,-2,1)$, and the equation of the plane is $2x+2y-z=0$ and $E_\Box(t)=6t^2+4t+1$. We know that $E_\Box(-1)=3$ represents the number of lattice points in the interior of the square. We notice that $3$ (and $E_\Box(-2)=17$)  was not an almost perfect number.  Hence we may want to  generalize that   sequence to {\em almost perfect squares in dimension $n$}.

What is interesting is that in the plane  $2x+2y-z=0$ we have another square: $\overline{u}=(2,-1,2)$ and $\overline{v}=(-1,2,2)$. Also, the square above can be written in terms of this one basically like in two dimensions: $u=\overline{u}-\overline{v}$ and $v=\overline{u}+\overline{v}$. For this square, $n=3(-2,-2,1)$ and so the Ehrhart polynomial is $E_{\{\overline{u},\overline{v}\}}(t)=3t^2+2t+1$. We see that $E_{\{\overline{u},\overline{v}\}}(-1)=2$ which is yet new for the almost perfect square sequence compared with dimension two.

If we take the twin vectors $u=(6,-2,3)$ and $v=(-2,3,6)$ we get $n=7(-3,-6,2)$ and so the Ehrhart polynomial is $E(t)=7t^2+2t+1$. There seems to be very few numbers that are not perfect squares in 3D: $7$, $14$, $23$,...

For an example of an irreducible square, for which $d$ and $d'$ in Theorem~\ref{epford3} are both greater than one, we refer to $u=5(8,12,9)$ and $v=17(0,-3,4)$. Its polynomial is $E(t)=(5t+1)(17t+1)$.

How does one construct such squares in $\mathbb R^3$? We observed in the proof of Theorem~\ref{epford3} that every such square is contained in a plane whose normal $n=(n_1,n_2,n_3)$ satisfies
\begin{equation}\label{characteristicequtionforsquares}
\ell^2=n_1^2+n_2^2+n_3^2.
\end{equation}

A solution of this type of Diphantine equation is usually referred in the literature as a {\em  Pythagorean quadruple}. For the number of primitive solutions of
(\ref{characteristicequtionforsquares}) in terms of $\ell$,  we refer the reader to a recent paper
of Werner $\rm H\ddot{u}rlimann$ (\cite{Hurlimann}) but also \cite{ch}. There are at least two other terms used for these quadruples of integers: cuboids (\cite{Hurlimann})
and {\it Lorenz quadruples} (see \cite{grahamLagariasMWC}). We include in the next tables all primitive solutions for odd values $\ell\in\{1,3,..19\}$:

$$\begin{tabular}{|c|c|}
  \hline
  $\ell$ & [a,b,c], $\gcd(a,b,c)=1$, $c$ \ even\\
   \hline
  1 & [1,0,0]\\
    3 & [1,2,2]\\
      5 & [3,4,0]\\
       7 & [3,6,2]\\
       9 & [1, 4, 8], [7, 4, 4] \\
  \hline
\end{tabular},
\begin{tabular}{|c|c|}
  \hline
  $\ell$ & [a,b,c], $\gcd(a,b,c)=1$ , $c$ even\\
   \hline
11& [9, 2, 6], [7, 6, 6]\\
13& [5, 12, 0], [3, 4, 12]\\
15& [5, 14, 2], [11, 2, 10]\\
17 & [15, 0, 8], [1, 12, 12], [9, 12, 8]\\
19& [1, 6, 18], [17, 6, 6], [15, 6, 10]\\
  \hline
\end{tabular}
$$

So, having a Pythagorean quadruple, how can we construct a related twin set of vectors? We include a partial answer to this question at this point, but we will include the complete solution later in this paper.

\begin{theorem}\label{constructionofsquares} Suppose that $u'=(a,b,c)$ satisfies $n_1a+n_2b+n_3c=0$, where $n_1^2+n_2^2+n_3^2=\ell^2$ with all variables involved being integers. Then there exist $v=(a',b',c')$ such that $\ell u'$ and $v$ define a lattice square in the plane of normal $n= (n_1,n_2,n_3)$.
\end{theorem}

\proof.  We define $v$ to be the cross-product of $u'$ and $n$. Clearly, $v$ is perpendicular to $n$ and so it is in the right plane.
Since $u'$ and $n$ are perpendicular,  $|v|=|u'||n|=\ell |u'|$.\eproof

 It is natural to look for the ``smallest" square in the plane $n_1a+n_2b+n_3c=0$. This leads us to the shortest vector problem (SVP), i.e., finding a non-zero vector in a lattice of minimum norm. An interesting problem, at this point, is to characterize all the values $\ell$ so that $\sqrt{\ell}$  appears as side-lengths for an lattice embedded square in $\mathbb R^m$, $m\ge 2$. If we denote this set by $\cal L$, from the two dimensions construction, we see that $\cal L$ is invariant under multiplication with numbers which are sums of two squares. For $m=3$, it is also clear  all the numbers of the form $4^k(8s-1)$ are not in  $\cal L$, since these are not representable as sums of three squares (Legendre's three-square theorem). Also,  $3$, $11$ or $19$ (see the sequence A223732 in OEIS) are not in $ \cal L$, since they have only one representation as sum of three squares
and these representations contain only  odd numbers. Hence the sum $aa'+bb'+cc'$ is also odd and so, it cannot be zero. We will describe the set $\cal L$ at the end of the subsection. We will see that in $\mathbb R^4$ and beyond, $\cal L=\mathbb N$.

We have the following parametrization for those squares whose side-lengths are natural numbers of the form $x^2+y^2+z^2+t^2$ (see \cite{ejinewparam} and \cite{spira}):

\begin{equation}\label{parametrization} u:=(2ty+2zx, 2tz-2yx, t^2-z^2-y^2+x^2), \
v:=(2zy-2tx, z^2-t^2+x^2-y^2, 2tz+2yx),
\end{equation}

\n having normal vector $n=(-x^2+t^2-y^2+z^2,-2(tx+zy),2(ty-zx))$, $|n|=x^2+y^2+z^2+t^2$.
 Next, we show that, as in \cite{ejinewparam} and \cite{spira}, the following similar parametrization for the solutions of (\ref{characteristicequtionforsquares}) takes place.
The proof is essentially the same as in \cite{spira} but we include it for the convenience of the reader.

\begin{theorem}\label{parametrizationforeq} (\cite{spira}) Every primitive solution of (\ref{characteristicequtionforsquares}), after a permutation of variables and change of signs, is given by $n_1=2(zy-tx)$, $n_2=2(tz+yx)$ and $n_3=z^2-t^2+x^2-y^2$, and $\ell =x^2+y^2+z^2+t^2$ for some integers $x$, $y$, $z$ and $t$.
\end{theorem}

\proof.\ As usual, a primitive solution of (\ref{characteristicequtionforsquares}) is one for which $\gcd(n_1,n_2,n_3)=1$. We must have $\ell$ odd since
otherwise $\gcd(n_1,n_2,n_3)\ge 2$. In this case, one of the $n_i$ must be odd and the other two even. Without loss of generality, let us assume that $n_1$ and $n_2$ are even. Then, we have  $$(\ell-n_3)(\ell+n_3)=n_1^2+n_2^2.$$  We know that $\ell+n_3=2\alpha$ and $\ell-n_3=2\beta$ are both even, and then the above equality becomes

$$\alpha \beta=(n_1/2)^2+(n_2/2)^2=A^2+B^2=(A+Bi)(A-Bi), $$

\n where $A=n_1/2$ and $B=n_2/2$.   A Gaussian prime  of the form $p=4k+3$ which divides $\alpha\beta$, must divide both $A$ and $B$. So, it cannot divide both $\alpha$ and $\beta$ because it then divides $\ell=(\alpha+\beta)/2$ and $n_3=(\alpha-\beta)/2$,
which implies that $\gcd(n_1,n_2,n_3)\ge p$. Hence, $\alpha$ and $\beta$ are both sums of two squares. If we let $\gcd(A+iB,\alpha)=t+iy$ and $\gcd(A+iB,\beta)=x+iz$. Let us make the observation that these equalities are defined up to a unit, i.e., $\pm 1$ or $\pm i$.

Since $\alpha$ is real,
$t-iy$ divides $\alpha$ and so $\alpha=(t^2+y^2)\alpha'$. But every prime factor of $\alpha$, appears as a factor of either $A+Bi$ or $A-Bi$. Taking into account the multiplicities we see that $\alpha$ divides $t^2+y^2$ and so $\alpha'=1$.  Similarly, we have $\beta=x^2+z^2$. Also, $t+iy$ divides $A+iB$ and also
$x+iz$ divides $A+iB$. Hence $A+iB=(t+iy)(x+iz)k$ and from here $A^2+B^2=(t^2+y^2)(x^2+z^2)|k|^2$ which together with what we have shown earlier forces $|k|=1$. So, by changing $x$ and $y$ we can assume that $k=1$. Hence we have $A=tx-zy$ and $B=xy+tz$. Then $\ell=\alpha+\beta=x^2+y^2+z^2+t^2$ and $n_3=\alpha-\beta=t^2+y^2-x^2-z^2$.\eproof

Theorem~\ref{parametrizationforeq} allows us to characterize $\cal L$. The result is not new, as we found recently, it is contained in \cite{gkms}. We include a proof of it, based on our development of the ideas.

\begin{theorem}\label{characterizationofl}  The set    of all $\ell$ so that $\sqrt{\ell}$  is the side-lengths for an  embedded square in $\mathbb Z^3$ is the set of
positive integers which are sums of two squares.
\end{theorem}

\proof.\  In one direction, i.e., $\cal L$ contains the set of all positive integers which are sums of two squares, we can observe that the lengths for squares in two dimensions
are of the form $\sqrt{a^2+b^2}$. For the other direction, if we start with an arbitrary square, $\cal S$,  in  $\mathbb Z^3$, let us suppose the side is $\sqrt{\ell}$.
As in the proof of  Theorem~\ref{constructionofsquares}, this square is contained in a plane of normal $n=(n_1,n_2,n_3)$, with $n_1^2+n_2^2+n_3^2=\ell^2$.
By Theorem~\ref{parametrizationforeq}, we can find $x$, $y$, $z$ and $t$, such that $n_1=2zy-2tx$, $n_2=2tz+2yx$ and $n_3=z^2-t^2+x^2-y^2$. Then  the parametrization
(\ref{parametrization}) gives a square in the same plane as $\cal S$, whose sides are $\ell$. We can then write the vector of this square in the basis
given by two of the  vectors in $\cal S$: $w=\alpha u+\beta v$ for some $\alpha$ and $\beta\in \mathbb Q$. Then, taking norms we get $\ell^2=\alpha^2 \ell+\beta^2 \ell$ so
$\ell=\alpha^2+\beta^2$. This shows that $\ell=A^2+B^2$ for some $A$, $B\in \mathbb Z$.\eproof

Theorem~\ref{constructionofsquares} and Theorem~\ref{characterizationofl} imply the next consequence which one can also prove elementary.

\begin{corollary}\label{pointsinaplane} Given a point $P=(x,y,z)\in \mathbb Z^3$ in the plane of equation $n_1x+n_2y+n_3z=0$ with
$n_1^2+n_2^2+n_3^2=\ell^2$ for some $\ell \in \mathbb N$, then the number $x^2+y^2+z^2$ is actually a sum of two squares.
\end{corollary}

\subsection{Squares in $\mathbb R^4$}\label{squares}

Perhaps, one of the simplest ways to construct squares in $\mathbb Z^4$ is to take two squares in two dimensions and ``add" them together, each on its own dimensions. In other words, the two vectors $u$ and $v$ that define the square, as we have seen before, are of the form $u=(a,b,c,d)$ and $v=(-b,a,-d,c)$. We have clearly $u\cdot v=0$ and  $\|u\|^2=\|v\|^2=a^2+b^2+c^2+d^2$. This simple example allows us to answer the question about the possible side-lengths of such a squares. By Lagrange's Four Square Theorem, we see that every natural number is a possible side-length of a square. Are these squares the most general situation that one can expect?  We observe that $u$ and $v$ are essentially the first two rows of the
following {\em  pseudo-orthogonal matrix} (every two rows are mutually orthogonal and have the same norm):

$$O_{a,b,c,d}=\left[
       \begin{array}{cccc}
         a & b & c & d \\
         -b & a & d & -c \\
        -c & -d & a & b \\
         -d & c & -b & a \\
       \end{array}
     \right].$$

It is not difficult to see that the product of matrices like these ($O^TO=tI$) is of the same form.  This allows one to define a certain multiplication on vectors in
$\mathbb R^4$ which is exactly the quaternion multiplication that we will be using latter.

{\bf Example 1.} Let us compute the Ehrhart polynomial for a particular case like this that which gives a new value for its sides. We have $u=(2,1,1,1)$ and $v=(-1,2,-1,1)$.
Using the same idea as in the three dimensions we are computing two normal vectors that define the orthogonal space of the plane generated by $u$ and $v$.
For a generic point $P=(x,y,z,t)\in \mathbb R^4$, the two equations, computed with the help of the cross-product in the three dimensions, are
$$3x-y-5z=0\ \text{and}\  2y+z-3t=0.$$

\n We observe that both $u$ and $v$ satisfy these equations. Since these two equations can be solved for $y$ and $t$, we can find two integer vectors which generate the lattice $\mathbb Z^4$ intersected with the plane defined by $u$ and $v$:  $y=3x-5z$ and $t=2x-3z$ which gives the two vectors $\alpha=(1,3,0,2)$ and $\beta=(0,-5,1,-3)$.
We know then that these two vectors form a fundamental domain, so we need to compute the area of the parallelogram generated by these two vectors: $A=|\alpha| |\beta|\sin(\gamma)$ where $\cos \gamma=\alpha\cdot \beta/|\alpha| |\beta|=-3/\sqrt{10}$. This gives $A=7$ and so, the  Ehrhart polynomial is $$E_{\Box}(t)=t^2+2t+1=(t+1)^2.$$

For a set of vectors $\cal S$ in $\mathbb R^4$ we denote as usual by ${\cal S}^{\bot}$ the set of all vectors $x\in \mathbb R^4$ perpendicular to every vector $v$
in $\cal S$, i.e., $${\cal S}^{\bot}=\{(x_1,x_2,x_3,x_4)\in \mathbb R^4|x_1v_1+x_2v_2+x_3v_3+x_4v_4=0\ \text{for all}\ v=(v_1,v_2,v_3,v_4),v\in \cal S\}.$$

\begin{theorem}\label{necessearycond} (i) Given  integer vectors $u=(u_1,u_2,u_3,u_4)$ and $v=(v_1,v_2,v_3,v_4)$ such that
$$0<\ell=u_1^2+u_2^2+u_3^2+u_4^2=v_1^2+v_2^2+v_3^2+v_4^2,\ \  u_1v_1+u_2v_2+u_3v_3+u_4v_4=0,$$

\n  then there exist an odd $k\in \mathbb N$ dividing $\ell$ and two vectors $w_1$ and $w_2$ with integer coordinates such that

$$w_1=(0,\alpha_1-\beta_1,\alpha_2-\beta_2,\alpha_3-\beta_3)\
\text{and} \ w_2=(\alpha_3-\beta_3, -\alpha_2-\beta_2,\alpha_1+\beta_1
,0)\ \ \text{with}$$

\begin{equation}\label{thetwoeq}
k^2=\alpha_1^2+\alpha_2^2+\alpha_3^2=\beta_1^2+\beta_2^2+\beta_3^2,
\end{equation}

\n   $\alpha_i$ and $\beta_i$ of the same parity and  $u,v\in
\{w_1,w_2\}^{\bot}$.
One can permute the coordinates of $u$ and $v$ and/or change their
signs in order to have $w_1$ and $w_2$ linearly independent.

(ii) If
$\gcd(\alpha_1,\beta_1,\alpha_2,\beta_2,\alpha_3,\beta_3)=1$ then
the volume of the fundamental domain of the minimal lattice
containing $u$ and $v$ is equal to $k$.

\end{theorem}

\n \proof. \ Using Lagrange's Identity, we get

$$\ell ^2=(\sum_{i=1}^4 u_i^2)(\sum_{i=1}^4  v_i^2)=(\sum_{i=1}^4u_iv_i)^2+\sum_{1\le i< j\le 4} (u_iv_j-u_jv_i)^2.$$

\n  Since $u_1v_1+u_2v_2+u_3v_3+u_4v_4=0$, this implies that

\begin{equation}\label{sixsquares}
\ell^2=\sum_{1\le i< j\le 4} (u_iv_j-u_jv_i)^2.
\end{equation}

If we denote by $\Delta_{ij}=(-1)^{i-j}(u_iv_j-u_jv_i)$,  it is not difficult to check that

\begin{equation}\label{planeeq}
\Delta_{12}\Delta_{34}-\Delta_{13}\Delta_{24}+\Delta_{14}\Delta_{23}=0.
\end{equation}

\n This helps us change (\ref{sixsquares}) into

\begin{equation}\label{threequares}
\ell ^2=(\Delta_{12}\pm \Delta_{34})^2+(\Delta_{13}\mp\Delta_{24})^2+(\Delta_{14}\pm\Delta_{23})^2.
\end{equation}

One can check that

\begin{equation}\label{planeeq2}
\begin{cases}u_1(0)+u_2\Delta_{34}+u_3\Delta_{24}+u_4\Delta_{23}=0\\
v_1(0)+v_2\Delta_{34}+v_3\Delta_{24}+v_4\Delta_{23}=0\\
u_1\Delta_{23}+u_2\Delta_{13}+u_3\Delta_{12}+u_4(0)=0\\
v_1\Delta_{23}+v_2\Delta_{13}+v_3\Delta_{12}+v_4(0)=0.
\end{cases}
\end{equation}

\n Let us observe that if $2$ divides $\ell$, then all $\Delta_{ij}$ are divisible by $2$, because all $u_i$ ($v_i$) are all even or all odd.
In order to prove the claim, we simplify (\ref{threequares}) by the greatest power of $2$ possible and take $\alpha_1=(\Delta_{12}+
\Delta_{34})/\ell'$, $\alpha_2=(-\Delta_{13}+\Delta_{24})/\ell'$,
$\alpha_3=(\Delta_{14}+\Delta_{23})/\ell'$, $\beta_1=(\Delta_{12}-
\Delta_{34})/\ell'$, $\beta_2=(-\Delta_{13}-\Delta_{24})/\ell'$, and
$\beta_3=(\Delta_{14}-\Delta_{23})/\ell'$, with $\ell'$ so that $\gcd(\alpha_1,\beta_1,\alpha_2,\beta_2,\alpha_3,\beta_3)=1$.

It is easy to check that if $\alpha_3\not = \beta_3$ then the two
vectors $w_1$ and $w_2$ are linearly independent and so the
orthogonal space  $\{w_1,w_2\}^{\bot}$ is two-dimensional. By
permuting coordinates and/or changing their signs, we can insure
that $\alpha_3-\beta_3=2\Delta_{23}/ \ell' \not=0$, since by
(\ref{sixsquares}) and the assumption that $\ell>0$, we see that
not all of the $\Delta_{ij}$ are equal to zero.

(ii) Let us denote by ${\cal L}$ the minimal lattice which
contains the given square $\{u,v\}$ and by $V$ its volume. As we observed in part (i)
this lattice is the same as $\mathbb Z^4\cap \{w_1,w_2\}^{\bot}$
if $\alpha_3\not = \beta_3$. Let us observe that the two vectors
$U:=(0,\Delta_{12},-\Delta_{13},\Delta_{14})$ and
$V:=(\Delta_{14},-\Delta_{24},\Delta_{34},0)$ are in ${\cal L}$.
This shows that
$2U/\ell'=(0,\alpha_1+\beta_1,\alpha_2+\beta_2,\alpha_3+\beta_3)$
and
$2V/\ell'=(\alpha_3+\beta_3,-\alpha_2+\beta_2,\alpha_1-\beta_1,0)$
are in ${\cal L}$ also. The area of the parallelogram determined by these two vectors
is given by the square root of the  Gramian determinant

$$\left|
    \begin{array}{cc}
    (\alpha_1+\beta_1)^2+(\alpha_2+\beta_2)^2+(\alpha_3+\beta_3)^2   & (\alpha_1+\beta_1)(\beta_2-\alpha_2)+ (\alpha_2+\beta_2)(\alpha_1-\beta_1)\\
   (\alpha_1+\beta_1)(\beta_2-\alpha_2)+ (\alpha_2+\beta_2)(\alpha_1-\beta_1)    & (\alpha_3+\beta_3)^2+(\alpha_2-\beta_2)^2+(\alpha_1-\beta_1)^2 \\
    \end{array}
  \right|=$$

$$=\left|
    \begin{array}{cc}
    2(k^2+\alpha_1\beta_1+\alpha_2\beta_2+\alpha_3\beta_3)   & 2(\alpha_1\beta_2-\beta_1\alpha_2)\\
2(\alpha_1\beta_2-\beta_1\alpha_2)       &   2(k^2+\alpha_3\beta_3-\alpha_2\beta_2-\alpha_1\beta_1) \\
    \end{array}
  \right|.$$

\n This means that the vectors  $U/\ell'$ and $V/\ell'$, still in $\cal L$, form a parallelogram of an area, which is the square root of

$$(1/4)[(k^2+\alpha_3\beta_3)^2-(\alpha_2\beta_2+\alpha_1\beta_1)^2-(\alpha_1\beta_2-\beta_1\alpha_2)^2]=$$

$$(1/4)[k^4+2k^2\alpha_3\beta_3+\alpha_3^2\beta_3^2-(k^2-\alpha_3^2)(k^2-\beta_3^2)]=k^2((\alpha_3+\beta_3)/2)^2.$$
It follows that $V$ divides $k|\alpha_3+\beta_3|/2$. In a similar way we can show that  $V$ divides $k|\alpha_i\pm \beta_i|/2$ for $i=1,2,3$.
Given the assumption that $\gcd(\alpha_1,\beta_1,\alpha_2,\beta_2,\alpha_3,\beta_3)=1$, we conclude that $V$ divides $k$.

In order to conclude that $V=k$, let us look at the  Gram determinant of the vectors $w_1/(\alpha_3-
\beta_3)$ and $w_2/(\alpha_3-
\beta_3)$, which if added to $\cal L$ extends the lattice and gives a basis for it. Since its volume can be determined by the same method as above, a similar calculation gives
$V'=k/|\alpha_3-\beta_3|$. This implies that $V/V'$ is an integer, or $k$ divides $V|\alpha_3-\beta_3|$. Since we can  obtain similarly that $k$ divides $V|\alpha_i\pm \beta_i|$,
we conclude that $V=k$. \eproof

{\bf Example 2:} As in the case of equilateral triangles (see \cite{ejieqtrinz4}), there is a converse of Theorem~\ref{necessearycond}. Let us take the first odd $\ell$ for which we get two essentially different representations as in  (\ref{thetwoeq}):
$$11^2=9^2+6^2+2^2=7^2+6^2+6^2.$$

We can then choose $\alpha_1=9$, $\beta_1=7$, $\alpha_2=6$, $\beta_2=6$, $\alpha_3=2$, and  $\beta_3=6$. Then the two vectors, defined in Theorem~\ref{necessearycond}, are $w_1=(0,2,0,-4)$ and $w_2=(-4,-12,16,0)$. Then the space $\{w_1,w_2\}^{\bot}$ is defined by the equations $y-2t=0$ and $x+3y-4z=0$. These can be simplified to $y=2t$ and $x=4z-6t$. Because the generic vector in  $\{w_1,w_2\}^{\bot}\cap \mathbb Z^4$ is $u=(4z-6t,2t,z,t)=z(4,0,1,0)+t(-6,2,0,1)$ we can calculate, as before the ``volume " of the fundamental domain and obtain indeed $11$.   Then we get a quadratic form $QF(z,t)=(4z-6t)^2+(2t)^2+z^2+t^2=17z^2-48zt+41t^2$ that should lead to the solutions we need solving the Diophantine equation $QF(z,t)=11k$. It turns out that the smallest multiple for which we have solutions is $k=11$ and then two of the vectors which define the square are $u=(-4,8,5,4)$ and $v=(10,2,4,1)$. Its Ehrhart polynomial is then $E_\Box(t)=11t^2+2t+1$. We observe that we get new terms, such as $10$, $94$, and $266$ for instance, for the sequence of almost perfect squares in dimension 4, compared to dimension $2$. This is an example when all the $\Delta_{ij}$ are divisible by $11$.

As we have seen in (\cite{ejieqtrinz4}), the following result gives a certain converse to Theorem~\ref{necessearycond}.

\begin{theorem}\label{maintheorem} Given $k$ odd, and two different representations
$$k^2=a^2+b^2+c^2=a'^2+b'^2+c'^2,\ \text{with}\ \gcd(a,b,c,a',b',c')=1,\ c'>c, $$

\n  and $a$, $a'$ both odd. Then if we set $\Delta_{12}=\frac{a'-a}{2}$, $\Delta_{34}=\frac{a+a'}{2}$, $\Delta_{13}=-\frac{b'-b}{2} $, $\Delta_{24}=\frac{b+b'}{2}$,
$\Delta_{14}=\frac{c+c'}{2}$, and $\Delta_{23}=\frac{c'-c}{2}$  the equations (\ref{planeeq}) and (\ref{threequares}) are satisfied.
Moreover, the two  dimensional space
$\cal S$ of all vectors $[u,v,w,t]\in \mathbb Z^4$, such that

\begin{equation}\label{planeeqsc}
\begin{cases}(0)u+\Delta_{34}v+\Delta_{24}w+\Delta_{23}t=0\\
 \Delta_{23}u+\Delta_{13}v+\Delta_{12}w+(0)t=0
\end{cases}
\end{equation}

\n contains a family of squares.
\end{theorem}

\n \proof.\  By construction we see that (\ref{threequares}) is true. Since $a$, $a'$ are odd, the others must be even and so, all numbers defined above are integers.  We have then indeed  $$\Delta_{12}\Delta_{34}-\Delta_{13}\Delta_{24}+\Delta_{14}\Delta_{23}=(1/4)[a'^2-a^2-(b^2-b'^2)+(c'^2-c^2)]
=0.$$

One can also check that $k^2=\sum_{i<j}\Delta_{ij}^2$.  We see that the assumption $\Delta_{23}>0$ insures that the equations (\ref{planeeqsc}) define a two dimensional space in $\mathbb R^4$.
We can solve the equations (\ref{planeeqsc}) for $t$ and $u$:

$$t=-\frac{\Delta_{34}v+\Delta_{24}w}{\Delta_{23}}, \ u= -\frac{\Delta_{13}v+\Delta_{12}w}{\Delta_{23}},\ \text{with}\ v,w\in \mathbb Z.$$

\n If we denote a generic point $P\in\mathbb R^4$ in the plane (\ref{planeeqsc}), in terms of $v$ and $w$, i.e.,

$$P(v,w)=[-\frac{\Delta_{13}v+\Delta_{12}w}{\Delta_{23}}, v,w, -\frac{\Delta_{34}v+\Delta_{24}w}{\Delta_{23}}],$$

\n we observe that for two pairs of the parameters, $(v,w)$ and $(v',w')$, we obtain $\Delta_{23}'=wv'-vw'$,
$$\Delta_{12}'=v(-\frac{\Delta_{13}v'+\Delta_{12}w'}{\Delta_{23}})+\frac{\Delta_{13}v+\Delta_{12}w}{\Delta_{23}}v'=\Delta_{12}\Delta_{23}'/\Delta_{23},$$
\n and similarly $\Delta_{13}'=\Delta_{13}\Delta_{23}'/\Delta_{23}$, etc.

We want show that a non-zero square the form $OP(v,w)P(v',w')Q$ ($OQ=OP(v,w)+OP(v',w')$) exists for some integer values of  $v$, $w$, $v'$ and $w'$.
We need to have $OP(v,w)^2=OP(v',w')^2=\ell$ and $P(v,w)\cdot P(v',w')=0$.

If such a square exists,   by Theorem~\ref{necessearycond}, we may want to  $\ell$ of the form  $\ell=k\ell'$ for some $\ell' \in \mathbb N$.

By Lagrange's identity, calculations similar to those  in the derivation of (\ref{sixsquares}), show that if
$|P(v,w)|^2=|P(v',w')|^2=k\ell'$ then

$$(P(v,w)\cdot P(v',w'))^2=k^2\ell'^2-(\sum_{i,j}\Delta_{ij}^2)\left(\frac{\Delta_{23}'}{\Delta_{23}}\right)^2=k^2\ell'^2-k^2\left(\frac{\Delta_{23}'}{\Delta_{23}}\right)^2.$$

\n Hence, to get $P(v,w)\cdot P(v',w')=0$, it is enough to have $\Delta_{23}'=\ell' \Delta_{23}$ (provided that we have already shown that $|P(v,w)|^2=|P(v',w')|^2=k\ell'$).
In other words, a square exists if there exist integer pairs $(v,w)$ and $(v',w')$ such that

\begin{equation}\label{equivfrmulation}
|OP(v,w)|^2=|OP(v',w')|^2=k\frac{\Delta_{23}'}{\Delta_{23}}.
\end{equation}

\n The essential object in this analysis is the quadratic form

$$QF(v,w):=|OP(u,v)|^2=\left(\frac{\Delta_{34}v+\Delta_{24}w}{\Delta_{23}}\right)^2+v^2+w^2+\left(\frac{\Delta_{13}v+\Delta_{12}w}{\Delta_{23}}\right)^2,$$

or

$$QF(v,w)=\frac{(\Delta_{34}^2+\Delta_{13}^2+\Delta_{23}^2)v^2+2(\Delta_{34}\Delta_{24}+\Delta_{13}\Delta_{12})vw
+(\Delta_{24}^2+\Delta_{12}^2+\Delta_{23}^2)w^2}{\Delta_{23}^2}.$$

\n We use Lagrange's identity again, which we will write in the form

$$(\alpha^2+\beta^2+\gamma^2)(\zeta^2+\eta^2+\theta^2)=(\alpha\zeta-\beta\eta+\gamma\theta)^2
+(\alpha\eta+\beta\zeta)^2+(\alpha\theta-\gamma\zeta)^2+(\beta\theta+\gamma\eta)^2.$$

\n Then the determinant of the form $QF(v,w)$ (excluding its denominator) is equal to

$$\begin{array}{c}
-\Delta/4=(\Delta_{34}^2+\Delta_{13}^2+\Delta_{23}^2)(\Delta_{12}^2+\Delta_{24}^2+\Delta_{23}^2)-(\Delta_{34}\Delta_{24}+\Delta_{13}\Delta_{12})^2=\\ \\
(\Delta_{12}\Delta_{34}-\Delta_{13}\Delta_{24}+\Delta_{23}^2)^2+(\Delta_{12}\Delta_{23}-\Delta_{34}\Delta_{23})^2+(\Delta_{13}\Delta_{23}+\Delta_{24}\Delta_{23})^2\Rightarrow\\ \\
-\Delta/4 =\Delta_{23}^2\left[(\Delta_{23}-\Delta_{14})^2+(\Delta_{12}-\Delta_{34})^2+(\Delta_{13}+\Delta_{24})^2 \right]=\Delta_{23}^2(k^2)\\ \\  \Rightarrow \Delta=-(2k\Delta_{23})^2.
\end{array}$$

\n This implies that for $v_0= -(\Delta_{34}\Delta_{24}+\Delta_{13}\Delta_{12})$ and $w_0=\Delta_{13}^2+\Delta_{23}^2+\Delta_{34}^2$ we get

$$QF(v_0,w_0)=k^2(\Delta_{34}^2+\Delta_{13}^2+\Delta_{23}^2),$$

\n and

\begin{equation}\label{qf}
QF(v,w)=\frac{(w_0v-v_0w)^2+k^2w^2\Delta_{23}^2}{\Delta_{23}^2w_0}.
\end{equation}

\n So, by (\ref{equivfrmulation}) and (\ref{qf}) we need to have solutions $(v,w)$ and $(v',w')$ of

 \begin{equation}\label{suffices}
(w_0v-v_0w)^2+k^2w^2\Delta_{23}^2 =(w_0v'-v_0w')^2+k^2w'^2\Delta_{23}^2 =k\Delta_{23}w_0(wv'-vw').
 \end{equation}

 \n In general, if the Diophantine quadratic equation $x^2+y^2=n$ has a solution $(x,y)$, then it has other solutions such as $(y,-x)$.

\n  Given $(v,w)$, let us assume that $(v',w')$ gives exactly this other solution ($x=w_0v-v_0w$, $y=kw\Delta_{23}$) when substituted in (\ref{suffices}):

$$\ds
\begin{cases}\ds w_0v'-v_0w'=kw\Delta_{23}\\ \\
\ds k\Delta_{23}w'=-(w_0v-v_0w).
\end{cases}
$$

\n This allows us to solve for $v'$ and $w'$ in order to calculate $\Delta_{23}'$:

$$v'=\frac{v_0^2w-v_0w_0v+k^2\Delta_{23}^2w}{2k\Delta_{23}w_0},\ w'=\frac{v_0w-w_0v}{k\Delta_{23}},$$

$$\Delta_{23}'=wv'-vw'=\frac{(w_0v-v_0w)^2+k^2w^2\Delta_{23}^2}{k\Delta_{23}w_0}.$$

\n This means that (\ref{suffices}) is automatically satisfied with this choice of $(v',w')$. In fact, we have shown that, for every integer values of
$(v,w)$, such that

\begin{equation}\label{criticaleq}
(w_0v-v_0w)^2+k^2w^2\Delta_{23}^2 =(kw_0\ell') \Delta_{23}^2
\end{equation}

\n taking $(v',w')$ as above, we automatically get  a square determined by vectors $OP(v,w), P(v',w')\in \mathbb Q^4$, which may not have integer coordinates.
In order for (\ref{criticaleq}) to have rational solutions, in $v$ and $w$, we need to have $kw_0\ell'$ a positive integer which is a product of
primes of the form $4s+1$, and all the other primes, in its prime factorization, have even exponents. Clearly, there exists a smallest positive integer $\ell'$ with this property. \eproof

\begin{figure}
\centering
\includegraphics[scale=.3]{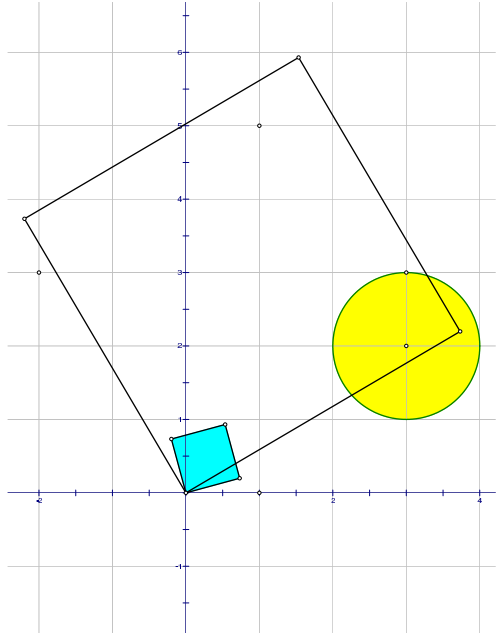}
\caption{Minimal Squares}
\label{fig:fig3}
\end{figure}

As we have seen in Proposition 2.2 (\cite{ejieqtrinz4}), we can similarly show the existence in every plane of equations (\ref{planeeqsc}),
of a square that is {\it minimal}, in the sense that the side-lengths are the smallest possible.   Then, a similar result to  Proposition 2.2 (\cite{ejieqtrinz4}) takes place: {\em every square in the same plane can be written in terms of a square that is minimal (same formulae as in 2 dimensions).}
For a proof we just invite the reader to study Figure~\ref{fig:fig3} (a proof ``without words" exercise).

As in dimension $3$, we also have parametric formulae for squares, which are more general than the ones at the beginning of the subsection.
The two vectors are

\begin{equation}\label{parmetrization4}
\begin{array}{c}
u=\pm (-ta-zb+yc-xd, -tb+za-yd-xc, -tc-zd-ya+xb, -td+zc+yb+xa)\\ \\
v=\pm (ax-by-cz-dt, ay+bx+ct-dz, az-bt+cx+dy, at+bz-cy+dx)
\end{array}
\end{equation}
\n and one can check that $u\cdot v=0$ and $|u|^2=|v|^2=(x^2+y^2+z^2+t^2)(a^2+b^2+c^2+d^2)$.
The two relations as in $(\ref{thetwoeq})$ can be simplified by $(x^2+y^2+z^2+t^2)^2$ and $(a^2+b^2+c^2+d^2) ^2$ respectively and reduced to
\begin{equation}\label{thetwoequationforparam}
\begin{array} {c}(x^2+y^2+z^2+t^2)^2=(x^2+y^2-z^2-t^2)^2+(2xt+2yz)^2+(2ty-2xz)^2,\\ \\   (a^2+b^2+c^2+d^2) ^2=(a^2+b^2-c^2-d^2)^2+(2ad-2bc)^2+(2ac+2bd)^2.
\end{array}
\end{equation}

These relations can be reduced even farther depending upon what $q_1$ and $q_2$ are. A simple observation here is that
the simplified relations (\ref{thetwoequationforparam}) depend only on the plane containing the square and every other
square contained in the same plane has the same simplified relations. This implies that the minimal square contained
in the plane has the sides at most $\sqrt{k_1k_2}$, where the $k_i^2=A_i^2+B_i^2+C_i^2$ are the primitive versions of (\ref{thetwoequationforparam}).

\vspace{0.1in}

{\bf Example 3.}  We start with two equalities as in (\ref{thetwoequationforparam}): $9^2=7^2+4^2+4^2$ and $15^2=11^2+10^2+2^2$.
This implies that $45^2=33^2+30^2+6^2=35^2+20^2+20^2$. As in Theorem~\ref{maintheorem}, we can take $a=33$, $a'=35$,
$b=30$, $b'=20$, $c=6$ and $c'=20$. Then $\Delta_{12}=1$, $\Delta_{34}=34$, $\Delta_{13}=5 $, $\Delta_{24}=25$,
$\Delta_{14}=13$, and $\Delta_{23}=7$, and the equations of the plane are $34v+25w+7t=0$ and \\
$7u+5v+w=0$. The system can be solved easily in terms of $u$ and $v$: $w=-7u-5v$ and $t=25u+13v$.
Then two vectors that generate the lattice in this plane are $(1,0,-7,25)$ and $(0,1,-5,13)$.
Then the ``volume" of the fundamental domain is $45$. The quadratic form is $Q(u,v)=675u^2+720uv+195v^2$
or $Q(u,v)=15(45u^2+48uv+13v^2)$. So, the equation $Q(u,v)=45\ell$ is equivalent to $(13v+24u)^2+9u^2=3(13)\ell$.
Then it is clear that $v$ has to be a multiple of $3$ and if we take $\ell=3$, $v=3$ and $u=-2$ we obtain a solution.
This gives ${\bf u}=(-2,3,-1,-11)$, and $v=6$ and $u=-3$ gives ${\bf v}=(-3,6,-9, 3)=3(-1,2,-3,1)$ (which is indeed a square). This means that the Ehrhart polynomial
of this square is $E_\Box(t)=3t^2+4t+1=(t+1)(3t+1)$, $t\in \mathbb N$. Taking the parameterizations as in Theorem~\ref{parametrizationforeq}
for the two equalities, $9^2=7^2+4^2+4^2$ and $15^2=11^2+10^2+2^2$, we observe that this square is covered by the parametrization
(\ref{parmetrization4}) by taking $a=x=z=1$, $b=d=t=2$, $c=0$, and $y=3$.  It is surprising that such a square with a big size for its side lengths
has no lattice points in the interior.

In order to understand better what is happening we will reformulate everything in terms of quaternions.

We remind the reader that the Hamilton quaternion algebra over the real numbers, denoted by ${\mathbb H}(\mathbb R)$,
is the associative unitary algebra given by the requirements:

(I) ${\mathbb H}(\mathbb R)$ is the free $\mathbb R$-module over the symbols $i$, $j$, and $k$, with $1$ the multiplicative unit;

(II) $i^2=j^2=k^2=-1$, $ij=-ji=k$, $jk=-kj=i$ and $ki=-ik=j$.

\vspace{0.1in}

\n If $q=x+yi+zj+tk\in {\mathbb H}(\mathbb R)$ the conjugate of $q$ is $\overline{q}=x-yi-zj-tk$ and the norm of $q$ is $N(q)=x^2+y^2+z^2+t^2$.
Some standard notation is then naturally appearing: $Re(q)=x$ and $Im(q)=yi+zj+tk$.

By ${\mathbb H}(\mathbb Z)$ we denote the subset of quaternions whose components are all integers. We imbed $\mathbb Z^4$ into ${\mathbb H}(\mathbb Z)$ in the natural way:
$(x,y,z,t)\hookrightarrow x+yi+zj+tk$.  Also, we will think of $\mathbb R^3$ imbedded in ${\mathbb H}(\mathbb R)$ in a, more or less, natural way $(y,z,t)\hookrightarrow yi+zj+tk$; in other words, $\mathbb R^3$ is  the hyperplane $Re(q)=0$.

It is known that this norm is multiplicative, i.e. $N(q_1q_2)=N(q_1)N(q_2)$, and $\overline{q_1q_2}=\overline{q_2}\ \overline{q_1}$.

\n For the important results about the arithmetic of ${\mathbb H}(\mathbb Z)$, we recommend the reader the
recent treatment in \cite{dsv}.  Using the same terminology as in \cite{dsv}, a quaternion $q$ is called {\it odd}, if $N(q)$ is an odd number.

The parametrization in Theorem~\ref{parametrizationforeq} is basically equivalent to
\begin{equation}\label{quaternions} n_1i+n_2j+n_3k=q(\epsilon)\overline{q}, \ \text{where}\ \ q=x+yi+zj+tk,\ \text{and}\ \epsilon \in\{i,j,k\}.
\end{equation}
\n In case $\epsilon=k$,  $n_1=2(xz+yt)$, $n_2=2(zt-xy)$, and $n_3=x^2-y^2-z^2+t^2$.

The parametrization (\ref{parmetrization4}) is equivalent to

\begin{equation}\label{quaternions2}u=q_1\epsilon_2 \overline{q_2}\ \ \text{and}\ \  v=q_1\epsilon_3 \overline{q_2}, \ \text{where}\ \ q_1=x+yi+zj+tk, q_2=a+bi+cj+dk,\end{equation}

\n with $\epsilon_2$ and $\epsilon_3$  fixed in $\{i,j,k\}$. In what follows we are going to take $\epsilon_2=j$ and $\epsilon_3=k$.
This square is in a plane defined by equations similar to (\ref{planeeqsc}) obtained from the two Pythagorean quadruples defined by $q_1$ and $q_2$
as in (\ref{quaternions}) using $\epsilon=\epsilon_1$ such that $\{\epsilon_1,\epsilon_2,\epsilon_3\}=\{i,j,k\}$.

\n {\bf Observation 1:} If we substitute $q=\widetilde{q}(\alpha + \beta k)$ into (\ref{quaternions}) for $\epsilon=k$, we see that because
$(\alpha  +\beta k )k=k (\alpha +\beta k) $, we have
$$n_1i+n_2j+n_3k=\widetilde{q} (k)(\alpha +\beta k)(\alpha -\beta k)\overline{\widetilde{q}}=(\alpha^2+\beta^2)\widetilde{q} (k)\overline{\widetilde{q}}=(\alpha^2+\beta^2)(\widetilde{n_1}i+\widetilde{n_2}j+\widetilde{n_3}k). $$
This implies that $n_i=\widetilde{n_i}(\alpha^2+\beta^2)$, $i=1,2,3$, which means that the primitive solutions of the equation (\ref{characteristicequtionforsquares})
have to arise from quaternions $q$ which have no right factors of the form $\alpha + \beta k$.
We obtain the same conclusion if $q=\widetilde{q}(1+i)$ or $q=\widetilde{q}(1+j)$, because $(1+i)k=-j(1+i)$ and $(1+j)k=i(1+j)$.  We have the following converse of this observation.

\begin{proposition}\label{converseparam} If in (\ref{quaternions}), we have $gcd(n_1,n_2,n_3)=n>1$,  i.e. if $n_1=2(xz+yt)$, $n_2=2(zt-xy)$ and $n_3=x^2-y^2-z^2+t^2$, are divisible by $n$, then $q=x+yi+zj+tk=(x'+y'i+z'j+t'k)\eta$ where $\eta\in \{1+ i,1+j,\alpha+k\beta\} $ for some integers $\alpha$ and $\beta$.
\end{proposition}

\ \proof.  Let us pick a prime $p$ which divides $n$.  Since $$n_1^2+n_2^2=4[(xz+yt)^2+(zt-xy)^2]=4(x^2+t^2)(y^2+z^2)$$ is divisible by $n^2$ it is divisible by $p$. First, we assume that $p>2$. Then $p$ divides either
$A:=x^2+t^2$ or $B:=y^2+z^2$. Since $n_3=A-B$ and $p$ divides $n_3$, we must have, in fact, $A$ and $B$ both divisible by $p$. If $p=2$, then $A$ and $B$ are either both even or both odd. If $A$ and $B$ are both even, then we have the same conclusion as in the case $p>2$. In the case $A$ and $B$ are both odd, we have, lets say
$x=2x'$, $y=2y'$, $z=2z'+1$, $t=2t'+1$. This implies that
$$q=x+yi+zj+tk=2x'+2y'i+2z'j+2t'k+j+k=(x'+y'i+z'j+t'k)(1-i)(1+i)+k(1+i)=\widetilde{q}(1+i)$$  which proves our claim.
Similar conclusion can be drawn if $x$ and $z$,  $t$ and $y$, or $t$ and $z$ are even.

If $p$ is a prime of the form $4k+3$, then automatically $x=px'$, $z=pz'$, $t=pt'$, and $y=py'$, which shows that
$q=x+yi+zj+tk=(x'+y'i+z'j+t'k)p$ and the conclusion of our proposition follows, with $\alpha=p$ and $\beta=0$.

Finally, if $p$ is a prime of the form $4k+1$, then $p=\alpha^2+\beta^2$. Because $A=(x+tk)(x-tk)$ and
$\eta=\alpha + \beta k$ is a Gaussian prime integer, it divides $x+tk$ or $x-tk$. Without loss of generality we may assume that $\eta$ divides $x+tk$.
If this is not the case we continue with $\overline{\eta}$. Next we observe that since $p$ divides $xz+yt$ and $zt-xy$, then
$p$ divides $(x+tk)(y+zk)=(xy-tz)+(ty+xz)k$. This implies that if $\overline{\eta}$ does not divides $x+tz$, then $\overline{\eta}$ must  divide $y+zk$.
Since $\eta$ divides $B=(y+zk)(y-zk)$, then
$\eta$ divides either $y+zk$ or $y-zk$. Let us assume first that $\eta$ divides $y-zk$. Then we can write

$$q=x+yi+zj+tk=x+tk+i(y-zk)=(x'+t'k)\eta+i(y'+z'k)\eta=\widetilde{q} \eta.$$

If $\eta$ divides $y+zk$, then $\overline{\eta}$ divides $y-zk$. If $\overline{\eta}$ divides $z+tz$ then we proceed as above, with a small change:

$$q=x+yi+zj+tk=x+tk+i(y-zk)=(x'+t'k)\overline{\eta}+i(y'+z'k)\overline{\eta}=\widetilde{q} \overline{\eta}.$$

 If $\overline{\eta}$  does not divides $x+tz$, then we have shown above that in this case, $\overline{\eta}$ must divide $y+zk$,  which is the same thing as $\eta$ dividing $y-zk$. This puts us in the same position as above. \eproof

\n {\bf Observation 2:} If we substitute $q_1=\widetilde{q_1}(\alpha + \beta i)$ into (\ref{quaternions2}), we see that because
$(\alpha  +\beta i )k=\alpha  k -\beta j $ and $(\alpha  +\beta i )j=\alpha  j +\beta k $, we have $\widetilde{u}=\alpha u - \beta v$  and $\widetilde{v}=\beta u + \alpha v$.
This is saying that the twin pair $\widetilde{u}$ and $\widetilde{v}$ is in the same plane as $u$ and $v$. A similar statement can be obtained if we substitute
$q_2=\widetilde{q_2}(\alpha + \beta i)$ into (\ref{quaternions2}).

Let us recall, for the convenience of the reader, the statement of Lemma 2.6.5 in \cite{dsv}, with the only difference that the important factors appear on the right
(so,  this new statement follows by conjugation and a symmetry type transformation $i\to -i$, $j\to -j$, and $k\to -k$).

\begin{lemma}\label{dsvlemma}  ({\bf \cite{dsv}})
Every integer quaternion $q\in {\mathbb H}(\mathbb Z)$ has unique factorization
\begin{equation}\label{factorizationofquaternions}
q=2^{\ell_q } \widetilde{q}\eta ,
\end{equation}
\n where $\widetilde{q}\in {\mathbb H}(\mathbb Z)$ is odd, $\eta_q\in \{1,1+i,1+j,1+k,(1+j)(1+i),(1-k)(1+i)\}$ and for some  non-negative $\ell_q\in \mathbb Z$.
\end{lemma}

This new statement follows by conjugation and a symmetry type transformation $i\to -i$, $j\to -j$, and $k\to -k$ from the original statement in \cite{dsv}.

\begin{theorem}\label{mostimportant} (i) We assume that $q_1$ and $q_2$ in the parametrization (\ref{quaternions2}) represented as in (\ref{factorizationofquaternions}), are not right-divisible by quaternions of the form $\pi=\alpha +\beta i$, $|\pi |>1$,  then the square in the parametrization (\ref{quaternions2}) is minimal.\par

\n (ii) The parametrization (\ref{quaternions2}) represents all the integer squares.
\end{theorem}

\ \proof. \ (i) We observe first that under our hypothesis, $\ell_{q_1}=\ell_{q_2}=0$ and  $\eta_{q_1}, \eta_{q_2}\in\{1+j,1+k\}$ in (\ref{factorizationofquaternions}). By way of contradiction, if the construction of the square $\cal S =\{u,v\}$ in (\ref{quaternions2}) is not minimal, there exists an integer square $\cal S_0 =\{u_0,v_0\}$ in the same plane in such a way $\cal S $ can be obtain from $\cal S_0 $ by a simple transformation which involves two integer parameters $\alpha$ and $\beta$:
$u=\alpha u_0+\beta v_0$ and $v=-\beta u_0+\alpha v_0$. We may assume that  $p=\alpha^2+\beta^2$ is a prime, otherwise we redefine the twin pair $\cal S_0 $ in such a way this condition is satisfied.
 Let us first assume that $p>2$. Since $q_1$ is not right-divisible by $\eta=\alpha +\beta i$ (which is a prime in ${\mathbb H}(\mathbb Z)$, $\eta$ odd quaternion), we conclude that $\gcd(q_1,\overline{\eta} )_r=1$ (the right greatest-common divisor as defined in \cite{dsv}). By Theorem 2.6.6 in \cite{dsv}, which is the equivalent of $\rm B\acute{e}zout$'s relation in ${\mathbb H}(\mathbb Z)$, we can find
$\gamma$ and $\delta$ in  ${\mathbb H}(\mathbb Z[\frac{1}{2}])$ such that

$$\gamma q_1+\delta \overline{\eta} =1.$$

If $p=2$, then  $|u|^2=|v|^2=N(q_1)N(q_2)=2|u_0|^2$ and so $2$ divides $N(q_1)$ or $N(q_2)$. Using Lemma~\ref{dsvlemma} we conclude that $q_1$
is of the form $q_1'(1+k)$ or  $q_1'(1+j)$ with $q_1'$ odd. So we can still obtain a  $\rm B\acute{e}zout$'s relation as above, with $q_1'$ and $\overline{\eta}$,
which can be easily transformed into one for $q_1$ and $\overline{\eta}$.

Solving for $u_0$, the above system, in terms of $u$ and $v$, we have $$u_0=(1/p)(\alpha u-\beta v)=(1/p) q_1( \alpha+i\beta) k \overline{q_2}.$$
Multiplying on the right with $(1/p) ( \alpha+i\beta) k \overline{q_2}$ the above $\rm B\acute{e}zout$'s relation, we obtain

$$\gamma u_0+\delta k \overline{q_2} =(1/p) ( \alpha+i\beta) k \overline{q_2}.$$

From this relation we see that $q_2=\widetilde{q_2}( \alpha-i\beta)$ for some $\widetilde{q_2}\in {\mathbb H}(\mathbb Z[\frac{1}{2}])$. Hence $2^kq_2=q_2'( \alpha-i\beta)$ for some $q_2'\in {\mathbb H}(\mathbb Z )$.  Using Lemma~\ref{dsvlemma}  again, we conclude that $q_2$ is right divisible by $\alpha-i\beta$ which is in contradiction with our hypothesis. Therefore, it
must be true that the $\cal S =\{u,v\}$ in (\ref{quaternions2}) is  minimal.

(ii) Let us start with an arbitrary square, $\cal S_0$. We then use
Theorem~\ref{necessearycond} to find the reduced equations of the plane that it is contained in, as those given by (\ref{planeeqsc}). Then we can apply the Theorem~\ref{parametrizationforeq} twice, and finally use the two sets of parameters (8 in total) to construct a square $\cal S_1$,
with formulae (\ref{parmetrization4}). The quaternions involved have the property in part (i) (both odd), and so the square $\cal S_1$ is minimal.
Then $\cal S_0$ can be written in terms of $\cal S_1$ using two numbers $\alpha$ and$\beta$.
Using the substitutions described before the statement the theorem, we can we can write $\cal S_0$ as in  (\ref{quaternions2}) by multiplying either of the quaternions
by  $\alpha+\beta i$.\eproof

\begin{theorem}\label{ehrhart4d} Assuming we have two odd integer quaternions $q_1$ and $q_2$ as in (ii) in Theorem~\ref{mostimportant}, then the fundamental domain of the integer lattice containing the square $\{u,v\}$ in (\ref{quaternions2}), has ``volume"  equal to $V=lcm (N(q_1)),N(q_2))$. The Ehrhart polynomial of the square  (\ref{quaternions2}) is
 $$E_\Box(t)=\gcd(N(q_1),N(q_2))t^2+(D_1+D_2)t+1,$$
 \n where $D_1=\gcd(u_1,u_2,u_3,u_4)$ and $D_2=\gcd(v_1,v_2,v_3,v_4)$.
\end{theorem}

\ \proof. \ This follows from Theorem~\ref{necessearycond} (part (ii)),  Proposition~\ref{converseparam} and the observation that $lcm(N(q_1),N(q_2))^2=a^2+b^2+c^2=a'^2+b'^2+c'^2$ for some integers $a$, $b$, $c$, $a'$, $b'$, and $c'$ with $\gcd(a,b,c,a',b',c')=1$.\eproof

\vspace{0.1in}

\n {\bf Observation 3:} To get back in dimension 3, we identified $\mathbb R^3$ with the subspace of quaternions $q$ for which  $Re(q)=0$. We have seen that
only one Pythagorean quadruple equation is necessary: $\Delta_{13}=\Delta_{12}=0$ which implies $\Delta_{14}=0$. To accomplish this,  we can take $q_1=q_2=q$ where $q$ is odd. This implies that the most general square in dimension 3 is of the form (\ref{parametrization}) combined with the
usual variations from dimension 2 in the respective plane:  $\widetilde{u}=\alpha u - \beta v$  and $\widetilde{v}=\beta u + \alpha v$ with $\{u,v\}$ as in (\ref{parametrization}).

There seems to be plenty of numerical evidence that the sequence of almost perfect squares in dimension four is the set of all non-negative integers. Let us denote the set of all such numbers by $\cal APS_4$. We include some partial result in this direction.

\begin{theorem}\label{almostperfectsquares} (i) The set $\cal APS_4$ contains all odd numbers.

(ii) The set $\cal APS_4$ contains all even numbers of the form $p-1$ with $p$ a prime $p\ge 11$.
\end{theorem}

\proof. \ (i) Let us take an odd $k$, $k\ge 1$, and a primitive Pythagorean quadruple representation $k^2=a^2+b^2+c^2$. Without loss of generality we may assume that
$a$ is odd and both $b$ and $c$ are even. If we multiply this equality by $2$,
we obtain $2k^2=2a^2+2b^2+2c^2=2c^2+(a+b)^2+(a-b)^2$. Hence we can build the square $u=(k,k,0,0)$ and $v=(c,-c,a+b,a-b)$. It is clear that $D_1=k$ and $D_2=1$ ($a-b$ and $a+b$ are odd)  as defined in Theorem~\ref{ehrhart4d}. Since the sides have length
$k\sqrt{2}$ and the two equation defining the plane containing $\{u,v\}$ are $k^2=(-c)^2+b^2+a^2=(-c)^2+a^2+b^2$, we conclude by Theorem~\ref{ehrhart4d} that the Ehrhart
polynomial is $E_{u,v}(t)=2kt^2+(k+1)t+1$. Hence $E_{u,v}(-1)=k$. This proves the first statement.

(ii) If $p$ is a prime $p\ge 11$, then we can find two distinct primitive Pythagorean quadruple representations $p^2=a^2+b^2+c^2=a'^2+b'^2+c'^2$ and construct a square as in (\ref{quaternions2}), whose Ehrhart polynomial is $E(t)=pt^2+(D_1+D_2)t+1$, where $D_i$ divide $p$. It is not possible to have $D_i=p$ since this puts us in the position where one of the vectors defining the square, say $u$, is of the form $u=(p,0,0,0)$. This leads essentially to only one equation which defines the plane containing $\{u,v\}$, but we have two different equations that define the plane.
Hence, since $E(-1)=p-1$ the claim in (ii) follows.\eproof

Numerical evidence shows that given an odd $k\ge 9$, one can find a square whose Ehrhart polynomial is $E(t)=kt^2+2t+1$. For instance, for $k=2015$ we found
$u=(-836, 584, -1592, -697)$ and $v=(-506, 1414, 203, 1328)$ that does the job.

\subsection{Cubes}

If a general lattice cube in $\mathbb R^3$ is given by the orthogonal matrix

\begin{equation}\label{orthogonalmatrix}
C_{\ell}=\frac{1}{\ell}\left[ \begin{array}{rrr}
          a_1 & b_1 & c_1 \\
          a_2 & b_2 & c_2 \\
          a_3 & b_3 & c_3\\
        \end{array}
      \right],
\end{equation}

\n with $a_i$, $b_i$ and $c_i$ integers satisfying $a_ia_j+b_ib_j+c_ic_j=\delta_{i,j}\ell^2$ for all $i$, $j$ in $\{1,2,3\}$), we define $d_i:=\gcd(a_i,b_i,c_i)$.
It is clear that the $d_i$ are divisors of $\ell$. Assuming that the matrix in (\ref{orthogonalmatrix}) is irreducible (no smaller $\ell$ can be used, which implies $\ell$ odd), then we have shown in \cite{ejilps} the following.

\begin{theorem}\label{ehrhartpforcubes}
Given a cube $C_{\ell}$ constructed from a matrix as in (\ref{orthogonalmatrix}), its Ehrhart polynomial is given by
\begin{equation}\label{ep}
L(C_{\ell},t)=\begin{array}{l} \ell^3 t^3+\ell(d_1+d_2+d_3)t^2+(d_1+d_2+d_3)t+1\ \text{or} \\ \\ (\ell t+1)[\ell^2 t^2+(d_1+d_2+d_3-\ell)t+1]
\end{array}, \ t\in \mathbb N.
\end{equation}
\end{theorem}

Similar  formula takes place for the cubes in dimension four. In general a square can be completed in various ways to a cube and even to a hypercube.
Let us assume that the cube is given by the first three vectors in the rows of the following matrix
\begin{equation}\label{orthogonalmatrix4d}
\left[ \begin{array}{rrrr}
          a_1 & b_1 & c_1 & d_1\\
          a_2 & b_2 & c_2 & d_2\\
          a_3 & b_3 & c_3& d_3\\
          a_4 & b_4 & c_4& d_4\\
        \end{array}
      \right],
\end{equation}
\n  with $a_i$, $b_i$, $c_i$ and $d_i$ integers $a_ia_j+b_ib_j+c_ic_j+d_id_j=\delta_{i,j}\ell $ for all $i$, $j$ in $\{1,2,3,4\}$). In this case, we define $D_i:=\gcd(a_i,b_i,c_i,d_i)$ and $\zeta_{ij}$ be the greatest common divisor of all the 2-by-2 determinants of the matrix
$$\left[ \begin{array}{rrrr}
          a_i & b_i & c_i & d_i \\
          a_j & b_j & c_j & d_j
        \end{array}
      \right],$$
\n for every $i$ and $j$ with $1\le i<j\le 3$.

\begin{theorem}\label{ehrhartpforcubes4d}
Given a cube $C_{\ell}$ constructed from a matrix as in (\ref{orthogonalmatrix4d}), its Ehrhart polynomial is given by
\begin{equation}\label{ep2}
L(C_{\ell},t)=\ell D_4 t^3+ \Delta t^2+\Delta' t+1,\ \
 t\in \mathbb N,
\end{equation}\n where $\Delta:=\zeta_{12}+\zeta_{13}+\zeta_{23}$ and $\Delta':=D_1+D_2+D_3$.
\end{theorem}

\n \proof. For the first coefficient the volume of the cube is $\ell \sqrt{\ell}$ and
 the volume of the fundamental domain of the lattice containing it is equal to$\sqrt{a_4^2+b_4^2+c_4^2+d_4^2}/D_4=\frac{\sqrt{\ell}}{D_4}$.
 For the second coefficient, the area of each of the six faces is $\ell$. Using the Theorem~\ref{mostimportant} to find the
 ``volume" of the fundamental domain of the lattice containing the faces determined by rows $i$ and $j$,  we obtain $\ell/\delta_{ij}$, which implies the stated value of $\Delta$.
Finally, for the coefficient $\Delta'$, we proceed like in the proof of Theorem~\ref{ehrhartpforcubes} (\cite{ejilps}). The main idea is essentially based on the fact that the topological equivalent in our cube to $[0,1)^3$ (containing $\kappa$  lattice points) has the fundamental domain property: the dilation with a factor $t$ contains exactly $t^3\kappa $ lattice points. This implies that $\kappa=\ell D_4$ and so

$$\ell D_4=1+(\ell D_4-\Delta+\Delta'-1)+\sum_{ij}(\zeta_{ij}-(D_i+D_j)+1)+\sum_{i}(D_{i}-1).$$
Solving for $\Delta'$ we obtain the stated value in the theorem. \eproof
\subsection{Hypercubes}

Hypercubes in dimension four can be constructed
using  quaternion techniques described earlier. A few orthogonal matrices in dimension $4$, together with the Ehrhart polynomial associated to the corresponding hypercube,  are included next:

$$ \underset{\begin{array}{c}  4t^4+8t^3+8t^2+4t+1\\ =(2t^2+2t+1)^2\end{array}}{\frac{1}{\sqrt{2}}\left[
  \begin{array}{cccc}
    1 & 1 & 0 & 0 \\
    1 & -1 & 0 & 0 \\
    0 & 0 & 1 & 1 \\
    0 & 0 & 0 & -1 \\
  \end{array}\right]}, \underset{\begin{array}{c} 9t^4+12t^3+6t^2+4t+1\\
 = (t+1)(3t+1)(3t^2+1)\end{array}}{\frac{1}{\sqrt{3}}\left[
  \begin{array}{cccc}
    1 & 1 & 1 & 0 \\
    -1 & 1 & 0 & 1 \\
    0 & -1 & 1 & 1 \\
    -1 & 0 & 1 & -1 \\
  \end{array}
\right]},
\underset{\begin{array}{c} 16t^4+16t^3+12t^2+4t+1\\
 = (1+2t+4t^2)^2\end{array}}{\frac{1}{2}\left[
  \begin{array}{cccc}
    1 & 1 & 1 & -1 \\
    -1 & 1 & 1 & 1 \\
    1 & -1 & 1 & 1 \\
    1 & 1 & -1 & 1 \\
  \end{array}
\right]},$$

{\small $$\underset{\begin{array}{c} 25t^4+20t^3+14t^2+4t+1\\
 = (1+2t+5t^2)^2\end{array}}{\frac{1}{\sqrt{5}}\left[
  \begin{array}{cccc}
    2 & 1 & 0 & 0 \\
    1 & -2 & 0 & 0 \\
    0 & 0 & 2 & 1 \\
    0 & 0 & 1 & -2 \\
  \end{array}
\right]}, \underset{36t^4+24t^3+8t^2+4t+1}{\frac{1}{\sqrt{6}}\left[
  \begin{array}{cccc}
    2 & 1 & 1 & 0 \\
    1 & -2 & 0 & 1 \\
    1 & 0 & -2 & -1 \\
    0 & 1 & -1 & 2 \\
  \end{array}\right]}, \underset{49t^4+28t^3+6t^2+4t+1}{\frac{1}{\sqrt{7}}\left[
  \begin{array}{cccc}
    2 & 1 & 1 & 1 \\
    1 & -2 & -1 & 1 \\
    1 & 1 & -2 & -1 \\
    1 & -1 & 1 & -2 \\
  \end{array}
\right]}, $$}

{\small $$\underset{\begin{array}{c} 81t^4+54t^3+18t^2+6t+1\\
 = (3t+1)^2(9t^2+1)\end{array}}{\frac{1}{3}\left[
  \begin{array}{cccc}
    3 & 0 & 0 & 0 \\
    0 & 2 & 2 & 1 \\
    0 & 2 & -1 & -2 \\
    0 & 1 & -2 & 2 \\
  \end{array}
\right]}, \underset{81t^4+36t^3+6t^2+4t+1}{\frac{1}{3}\left[
  \begin{array}{cccc}
    2 & 2 & 1 &  0\\
    2 & -2 & 0 &  1\\
    1 & 0 & -2 &  -2\\
    0 & 1 & -2 &  2\\
  \end{array}\right]},
   \underset{100t^4+40t^3+16t^2+4t+1}{\frac{1}{\sqrt{10}}\left[
  \begin{array}{cccc}
    2 & 2  & 1 & 1 \\
    2 & -2  & -1 & 1 \\
    1 & 1  & -2 & -2 \\
    1 & -1  & 2 & -2 \\
  \end{array}
\right]}, $$}

{\small $$\underset{\begin{array}{c} 100t^4+40t^3+24t^2+6t+1\\
 = (10t^2+2t+1)^2\end{array}}{\frac{1}{\sqrt{10}}\left[
  \begin{array}{cccc}
    3 & -1 & 0 & 0 \\
    1 & 3 & 0 & 0 \\
    0 & 0 & 3 & 1 \\
    0 & 0 & 1 & -3 \\
  \end{array}
\right]},  \underset{121t^4+44t^3+6t^2+4t+1}{\frac{1}{\sqrt{11}}\left[
  \begin{array}{cccc}
    3 & 1 & 1 & 0 \\
    1 & -3 & 0 & 1 \\
    1 & 0 & -3 & -1 \\
    0 & 1 & -1 & 3 \\
  \end{array}\right]}, \ \text{and}\ \underset{169t^4+53t^3+6t^2+4t+1}{\frac{1}{\sqrt{13}}\left[
  \begin{array}{cccc}
    2 & 2 & 2 & 1 \\
    2 & -2 & 1 & -2 \\
    2 & -1 & -2 & 2 \\
    1 & 2 & -2 & -2 \\
  \end{array}
\right]}.$$}

We see that several of these examples are cross polytopes (squares in two dimensions) and indeed the resulting polynomial is the
product of the smaller degree polynomials involved. Let us assume that in general we have

\begin{equation}\label{orthogonalmatrix4dn}
H_{\ell}=\frac{1}{\sqrt{\ell}}\left[ \begin{array}{rrrr}
          a_1 & b_1 & c_1 & d_1\\
          a_2 & b_2 & c_2 & d_2\\
          a_3 & b_3 & c_3& d_3\\
          a_4 & b_4 & c_4& d_4\\
        \end{array}
      \right],
\end{equation}

\n where $a_ia_j+b_ib_j+c_ic_j+d_id_j=\delta_{ij}\ell$ for $i,j\in\{1,2,3,4\}$. We define $D_i=\gcd(a_i,b_i,c_i,d_i)$ for $i\in\{1,2,3,4\}$ and assume that
$H_{\ell}$ is irreducible, i.e., $\gcd(D_1,D_2,D_3,D_4)=1$. The Ehrhart polynomial associated to the hypercube constructed in a natural way using the orthogonal matrix
 $H_{\ell}$ is then

 \begin{equation}\label{ehrhartforhypercube}
E_{H(\ell)}(t)=\ell^2 t^4+\alpha_1t^3+\alpha_2t^2+\alpha_3t+1.
 \end{equation}
\begin{proposition}\label{coeffc1} With the notation introduced earlier $\alpha_1=\ell(D_1+D_2+D_3+D_4)$.
\end{proposition}
\n \proof. Let us look at a 3-dimensional face, say generated by first three rows of (\ref{orthogonalmatrix4d}). Its volume is $\ell \sqrt{\ell}$ and
 the volume of the lattice containing it is equal to$\sqrt{a_4^2+b_4^2+c_4^2+d_4^2}/D_4=\frac{\sqrt{\ell}}{D_4}$. Since there are eight such faces the general theory gives the value claimed.\eproof

\begin{theorem}\label{thelasttheorem}
We have

\begin{equation}\label{thelastidentity}
\alpha_2=\alpha_3+\delta-\Delta,
\end{equation}
\n where $\delta=\sum_{1\le i<j\le 4}\zeta_{ij}$ and $\Delta=D_1+D_2+D_3+D_4$.
\end{theorem}

\n \proof.   We are going to use the same property of the part of hypercube that is topologically equivalent to $[0,1)^4$, as in  the proof of Theorem~\ref{ehrhartpforcubes4d}. The balancing of points in each of the corresponding sets of the partition

$$[0,1)^4=(0,0,0,0)\cup  (0,1)^4  \underset{\ all\ four\ 3D\ cubes} {\bigcup} (0,1)^3 \underset{\ all\ six\ 2D\ faces} {\bigcup} (0,1)^2\underset{\ all\ four\ 1D\ sides} {\bigcup} (0,1)   $$

\n gives (using Theorem~\ref{ehrhartpforcubes4d} )

$$\ell^2=1+(\ell^2-\alpha_1+\alpha_2-\alpha_3+1)+(\ell \Delta-2\delta+4\Delta-4)+\sum_{1\le i<j\le 4}(\zeta_{ij}-(D_i+D_j)+1)+\Delta-4. $$

\n From this, one can easily derive (\ref{thelastidentity}).\eproof

Looking at the examples we have so far, we notice that $\alpha_3=\Delta$, and so the Ehrhart polynomial of an hypercube would result into the follows simple form

\begin{equation}\label{last}E_{H(\ell)}(t)=\ell^2 t^4+\ell \Delta t^3+\delta t^2+\Delta t+1.
\end{equation}

This is indeed the case since it follows from Theorem 9.9 in \cite{beckAndRobins2015textbook}.
Of course, one may want to generalize all these concepts to dimensions bigger than four.

\end{document}